\documentclass[reqno,12pt]{amsart} 
\usepackage{amsmath}
\usepackage{amssymb}
\usepackage{amsthm}
\usepackage{graphicx, xcolor}
\usepackage[active]{srcltx} 
\newcommand\NoBlackBoxes{\global\overfullrule0pt}
\NoBlackBoxes
\theoremstyle{plain} 

\def\4{\kern1pt}

\def\6{\vphantom0}

\def\8{\kern-10pt}
\def\7#1{_{(#1)}}

\def\ffrac#1#2{\raise.5pt\hbox{\small$\4\displaystyle\frac{\,#1\,}{\,#2\,}\4$}}

\newcommand{\tc}{\textcolor{red}}
\begin{document}

\makeatletter
\let\serieslogo@\relax
\let\@setcopyright\relax

\def\speciallabelmark#1{\def\@currentlabel{#1}}
\makeatother


\title{Berry-Esseen Bounds\\ 
FOR TYPICAL WEIGHTED SUMS}

\author{S. G. Bobkov$^{1,4}$}
\thanks{1) School of Mathematics, University of Minnesota, USA;
Email: bobkov@math.umn.edu}
\address
{Sergey G. Bobkov \newline
School of Mathematics, University of Minnesota  \newline 
127 Vincent Hall, 206 Church St. S.E., Minneapolis, MN 55455 USA
\smallskip}
\email {bobkov@math.umn.edu} 
\thanks{\hskip4mm Research was partially supported by NSF grant DMS-1612961}

\author{G. P. Chistyakov$^{2,4}$}
\thanks{2) Faculty of Mathematics, University of Bielefeld, Germany;
Email: chistyak@math.uni-bielefeld.de}
\address
{Gennadiy P. Chistyakov\newline
Fakult\"at f\"ur Mathematik, Universit\"at Bielefeld\newline
Postfach 100131, 33501 Bielefeld, Germany}
\email {chistyak@math.uni-bielefeld.de}

\author{F. G\"otze$^{3,4}$}
\thanks{3) Faculty of Mathematics, University of Bielefeld, Germany;
Email: goetze@math.uni-bielefeld.de}
\address
{Friedrich G\"otze\newline
Fakult\"at f\"ur Mathematik, Universit\"at Bielefeld\newline
Postfach 100131, 33501 Bielefeld, Germany}
\email {goetze@mathematik.uni-bielefeld.de}

\thanks{4) Research was partially supported by the SFB 1283 at Bielefeld University}



\begin{abstract}
Under correlation-type conditions, we derive upper bounds of order 
$\frac{1}{\sqrt{n}}$ for the Kolmogorov distance between the distributions 
of weighted sums of dependent summands and the normal law.
\end{abstract}

\maketitle
\markboth{S. G. Bobkov, G. P. Chistyakov and F. G\"otze}{Berry-Esseen Bounds}

\def\theequation{\thesection.\arabic{equation}}
\def\E{{\mathbb E}}
\def\R{{\mathbb R}}
\def\C{{\bf C}}
\def\P{{\mathbb P}}
\def\H{{\rm H}}
\def\Im{{\rm Im}}
\def\Tr{{\rm Tr}}

\def\k{{\kappa}}
\def\M{{\cal M}}
\def\Var{{\rm Var}}
\def\Ent{{\rm Ent}}
\def\O{{\rm Osc}_\mu}

\def\ep{\varepsilon}
\def\phi{\varphi}
\def\F{{\cal F}}
\def\L{{\cal L}}

\def\be{\begin{equation}}
\def\en{\end{equation}}
\def\bee{\begin{eqnarray*}}
\def\ene{\end{eqnarray*}}


\vskip10mm
\section{{\bf Introduction}}
\setcounter{equation}{0}

\vskip2mm
\noindent
Given a random vector $X = (X_1,\dots,X_n)$ in $\R^n$ ($n \geq 2$), we consider 
the weighted sums
$$
S_\theta = \theta_1 X_1 + \dots + \theta_n X_n, \qquad 
\theta = (\theta_1,\dots,\theta_n) \in S^{n-1},
$$
parameterized by points of the unit sphere
$
S^{n-1} = \{\theta \in \R^n: \theta_1^2 + \dots + \theta_n^2 = 1\}.
$
According to the celebrated result by Sudakov [S], if $n$ is large, and 
if the covariance matrix of $X$ has bounded spectral radius,
the distribution functions $F_\theta(x) = \P\{S_\theta \leq x\}$
concentrate around a certain typical distribution function
given by the mean
\be
F(x) = \E_\theta F_\theta(x) \equiv \int_{S^{n-1}} F_\theta(x)\,d\mu_{n-1}(\theta), 
\qquad x \in \R,
\en
over the uniform probability measure $\mu_{n-1}$ on $S^{n-1}$. Although 
this theorem has a rather universal range of applicability (in contrast 
to the classical scheme of independent summands), the problem of possible 
rates of concentration, including the rates for the $\mu_{n-1}$-mean of the
Kolmogorov distance
$$
\rho(F_\theta,F) = \sup_x |F_\theta(x) - F(x)|,
$$
is rather delicate, and the answers depend upon correlation-type characteristics 
of the distribution of $X$. A natural characteristic is for example 
the maximal $L^p$-norm
$$
M_p = \sup_\theta \big(|\E S_\theta|^p\big)^{1/p}, \qquad p \geq 1.
$$

Moreover, if we want to study the approximation for most of $F_\theta$'s
by the standard normal distribution function
$$
\Phi(x) = \frac{1}{\sqrt{2\pi}} \int_{-\infty}^x e^{-y^2/2}\,dy, \qquad x \in \R,
$$
one is led to study another concentration problem -- namely rates for 
the distance $\rho(F,\Phi)$. To this aim, let us rewrite the definition (1.1) 
as $F(x) = \P\{r Z_n \leq x\}$ with
$$
r^2 = \frac{|X|^2}{n} = \frac{X_1^2 + \dots + X_n^2}{n} \qquad (r \geq 0),
$$
where the random variable $Z_n$ is independent of $r$ and has the same 
distribution as $\sqrt{n} \theta_1$ under $\mu_{n-1}$. Since $Z_n$ is 
close to being standard normal, $F$ itself is approximately normal, 
if and only if $r^2$ is nearly a constant, which translates into a weak law 
of large numbers for the sequence $X_k^2$. This property -- that 
the distribution of $r^2$ is concentrated around a point --
may be quantified by the variance-type functionals
$$
\sigma_{2p} = \sqrt{n}\, \big(\E\, |r^2 - 1|^p\big)^{1/p},
$$
which are expected to be of order 1 in reasonable situations (at least, 
they are finite, as long as $M_{2p} < \infty$). For example, if $|X|^2 = n$ a.s.,
we have $\sigma_{2p} = 0$. If the components $X_k$ are pairwise independent,
identically distributed, and with $\E X_1^2 = 1$, then
$$
\sigma_4^2 = \frac{1}{n}\,\Var(|X|^2) = \Var(X_1^2).
$$

It turns out that control of the two functionals, $M_3$ and $\sigma_3$ 
is sufficient to guarantee a Berry-Esseen type rate of normal approximation 
for $F_\theta$ on average, in analogy with
the Berry-Esseen theorem for independent identically distributed random variables.
Since the second moment for the typical distribution $F$ is equal to
$\E r^2$, a normalization condition for this moment is desirable.

\vskip5mm
{\bf Theorem 1.1.} {\sl If $\E\,|X|^2 = n$, then with some absolute constant $c$
\be
\E_\theta\, \rho(F_\theta,\Phi) \, \leq \, 
c\,(M_3^3 + \sigma_3^{3/2})\,\frac{1}{\sqrt{n}}.
\en
}

In the case of non-correlated random variables $X_k$'s, with mean zero and 
variance one, all $S_\theta$ have also mean zero and variance one, so that 
$M_2 = 1$. In many interesting examples, $M_3$ is known to be of the same order 
as $M_2$ (in particular, when Khinchine-type inequalities are available for 
linear functionals of $X$). In some other examples, the magnitude of $M_3$ 
is however much larger, and here control via $M_2$
is preferable, as the following assertion shows.

\vskip5mm
{\bf Theorem.1.2.} {\sl If $\E\, |X|^2 = n$, then for some absolute constant $c$
\be
\E_\theta\, \rho(F_\theta,\Phi) \, \leq \, 
c\,(M_2^2 + \sigma_2)\,\frac{\log n}{\sqrt{n}}.
\en
}

Thus, modulo an additional logarithmic factor, a Berry-Esseen type rate holds 
for this average under a second moment assumption, only.

For an illustration, consider the trigonometric system $X = (X_1,\dots,X_n)$ 
with components
\bee
X_{2k-1}(\omega) = \sqrt{2}\,\cos(k\omega),  & \\
X_{2k}(\omega) = \sqrt{2}\,\sin(k\omega), & \qquad -\pi < \omega < \pi, \ \ 
k = 1,\dots,n/2,
\ene
assuming that $n$ is even. They may be treated as random variables on the 
probability space $\Omega = (-\pi,\pi)$ equipped with the normalized Lebesgue 
measure $\P$, such that the linear forms
$$
S_\theta =  \sqrt{2}\, \sum_{k=1}^{n/2} 
\big(\theta_{2k-1} \cos(k\omega) + \theta_{2k} \sin(k\omega)\big)
$$
represent trigonometric polynomials of degree at most $\frac{n}{2}$. 
The normalization $\sqrt{2}$ is chosen for convenience only, since then $X$ is 
isotropic, so that $M_2 = 1$. Since also $\sigma_2 = 0$, by Theorem 1.2, 
most of the distributions $F_\theta$ of $S_\theta$ 
are approximately standard normal, and we have an upper bound
\be
\E_\theta\, \rho(F_\theta,\Phi) \, \leq \, c\,\frac{\log n}{\sqrt{n}}.
\en

The study of asymptotic normality for trigonometric polynomials has 
a long history, starting with results on lacunary systems due to Kac [Ka], 
Salem and Zygmund [S-Z1-2], Gaposhkin [G]; see also [B-G1-2], [A-B], [F], [A-E]. 
As we see, normality with an almost Berry-Esseen type rate 
remains valid for most choices of coefficients even without an assumption 
of lacunarity. One can show that the inequality 
(1.4) still holds for many other functional orthogonal systems as well, including, 
for instance, Chebyshev's polynomials on the interval $\Omega = (-1,1)$, 
the Walsh system on the Boolean cube $\{-1,1\}^n$. It holds as well for 
any system of functions of the form 
$X_k(\omega_1,\omega_2) = f(k\omega_1 + \omega_2)$,
$\omega_1,\omega_2 \in (0,1)$, where $f$ is 1-periodic and belongs to $L^4(0,1)$
(this is a strictly stationary sequence of pairwise independent random variables).
A common feature of all listed examples is that (1.4) may actually be reversed modulo 
a logarithmic factor, in the sense
that
$$
\E_\theta\, \rho(F_\theta,\Phi) \, \geq \, c\,\frac{1}{\sqrt{n}\,(\log n)^s}
$$
with some $s>0$. (However, we do not derive lower bounds here referring 
the interested reader to [B-C-G2]).

The conditions of Theorem 1.2 may be further relaxed in order to eliminate 
dependence on $\sigma_2^2$. This can be achieved by replacing it by the 
requirement of small probabilities for $\P\{|X-Y|^2 \leq n/4\}$, where $Y$ 
is an independent copy of $X$, cf. Theorem 6.3 below.
This extends the applicability of our results to further groups of examples, 
while replacing $\Phi$ by a certain mixture of centered Gaussian measures. 
More precisely, define $G$ to be the law of $rZ$, where $Z \sim N(0,1)$ 
is independent of $r = \frac{1}{n}\,|X|$. In particular, we have:

\vskip5mm
{\bf Theorem 1.3.} {\sl If the components $X_k$ of the random vector $X$ 
in $\R^n$ are independent, identically distributed, have mean zero and 
finite second moment, then
$$
\E_\theta\, \rho(F_\theta,G) \, \leq \, c\,\sqrt{\frac{\log n}{n}},
$$
where the constant $c$ depends on the distribution of $X_1$ only.
}

\vskip5mm
At first sight it seems surprising that an approximate Berry-Esseen type 
rate holds under no additional assumption beyond the finiteness of the 
second moment. Indeed, in the classical situation of equal coefficients, and when
$\E X_1 = 0$, $\E X_1^2 = 1$, the distributions $F_n$ of the normalized sums
$S_n = (X_1 + \dots + X_n)/\sqrt{n}$ may approach the standard normal law 
at an arbitrary slow rate: For any sequence $\ep_n \rightarrow 0^+$, one may
choose the distribution of $X_1$ such that 
$$
\rho(F_n,\Phi) \geq \ep_n
$$
for all $n$ large enough (cf. [M]). This shows that for typical coefficients, 
the distributions $F_\theta$ behave in a more stable way in comparison to 
$F_n$. This interesting phenomenon has been studied before.
For example, Klartag and Sodin [K-S] have shown in the i.i.d. case and under 
the 4-th moment assumption, that
$$
\E_\theta\, \rho(F_\theta,\Phi) \, \leq \, c\,\frac{\beta_4}{n}, \qquad 
\beta_4 = \E X_1^4,
$$
thus essentially improving the standard rate in the Berry-Esseen theorem
(see also [Kl]). 

The paper is organized as follows.
We start with comments on general properties of the moment and variance-type
functionals. Then we turn to the normal approximation for distributions of the
first coordinate on the sphere (with rate of order $1/n$), which is used in 
Section 4 to describe proper bounds on the distance from the typical distributions 
to the standard normal law. Proofs of both Theorems 1.1 and 1.2 rely upon 
the spherical Poincar\'e inequality and Berry-Esseen-type estimates in terms 
of characteristic functions.
The characteristic functions of the weighted sums are discussed separately 
in Section 5. Their properties are used in Section 6 to complete the proof 
of Theorem 1.2 (in a more general form). Theorem 1.1 is proved in Section 8, 
and in the last section we add some remarks concerning Theorem 1.3.


\vskip10mm
\section{{\bf Moment and variance-type functionals}}
\setcounter{equation}{0}

\vskip2mm
\noindent
First let us describe some basic properties of the functionals 
$M_p = M_p(X)$ and $\sigma_{2p} = \sigma_{2p}(X)$. We shall as well 
introduce a few additional functionals. Define
\be
m_p = m_p(X) = \frac{1}{\sqrt{n}}\, \Big(\E\, |\left<X,Y\right>|^p\Big)^{1/p}, \qquad
p \geq 1,
\en
where $Y$ is an independent copy of $X$. 

All these quantities do not depend on the
systems of coordinates, that is, $m_p(UX) = m_p(X)$ and $M_p(UX) = M_p(X)$ 
for any orthogonal linear map $U:\R^n \rightarrow \R^n$. 

We call $M_p$ the $p$-th moment of $X$. In case $M_2$ is finite, one may consider
the covariance operator (matrix) of $X$ which is defined by the equality
$$
\E \left<X,a\right>^2 = \left<Ra,a\right>, \qquad a \in \R^n.
$$
It is symmetric, positive definite, and has non-negative eigenvalues $\lambda_i$
($1 \leq i \leq n$). Choosing a system of coordinates such that $R$ is diagonal,
with  entries $\lambda_i$, we see that
\be
M_2^2 = \max_i \lambda_i, \qquad 
m_2^2 = \frac{1}{n}\, \E\left<X,Y\right>^2 =
\frac{1}{n}\, \sum_{i=1}^n \lambda_i^2,
\qquad \E\, |X|^2 = \sum_{i=1}^n \lambda_i.
\en

The random vector $X$ is called isotropic (or having an isotropic distribution), 
if the covariance matrix of $X$ is an identity, i.e.,
$$
\E \left<X,a\right>^2 = |a|^2, \qquad {\rm for \ all} \ \ a \in \R^n.
$$
In this case, $m_2 = M_2 = 1$, and $\E\,|X|^2 = n$.
Isotropic distributions are invariant under orthogonal transformations of the space.
Applying Cauchy's inequality, from (2.2) we immediately obtain:

\vskip5mm
{\bf Proposition 2.1.} {\sl For any random vector $X$ in $\R^n$ with $\E\,|X|^2 = n$,
we have $m_2 \geq 1$, where equality is attained, if and only if $X$ is isotropic.
}

\vskip5mm
The $p$-th moments of $X$ may easily be related to the moments of $|X|$.

\vskip5mm
{\bf Proposition 2.2.} {\sl Given $p \geq 2$, for any random vector $X$ in 
$\R^n$,
$$
(\E\, |X|^p)^{1/p} \leq M_p \sqrt{n}.
$$
}

If $X$ is isotropic, there is an opposite inequality
$
(\E\, |X|^p)^{1/p} \geq (\E\, |X|^2)^{1/2} = \sqrt{n}.
$

\vskip5mm
{\bf Proof.} By the rotational invariance of the uniform distribution on $S^{n-1}$,
we have 
$$
\E_\theta\,|\left<\theta,a\right>|^p = |a|^p\, \E_\theta\,|\theta_1|^p,
\qquad a \in \R^n,
$$
where $\E_\theta$ denotes the integral over the uniform measure $\nu_{n-1}$ .
Inserting here $a = X$, we get
$$
|X|^p\, \E_\theta\,|\theta_1|^p = \E_\theta\,|\left<X,\theta\right>|^p.
$$
Next, take the expectation with respect to $X$ and use
$\E\,|\left<X,\theta\right>|^p \leq M_p$ to arrive at the upper bound
$$
\E\,|X|^p \leq \frac{M_p^p(X)}{\E_\theta\,|\theta_1|^p}.
$$
Here, since $\E_\theta\, \theta_1^2 = \frac{1}{n}$, we have
$$
(\E_\theta\, |\theta_1|^p)^{1/p} \geq (\E_\theta\, |\theta_1|^2)^{1/2} 
= \frac{1}{\sqrt{n}}.
$$
\qed

\vskip2mm
{\bf Corollary 2.3.} {\sl $m_p \leq M_p^2$ for any $p \geq 2$.
}

\vskip5mm
Indeed, let $Y$ be an independent copy of the random vector $X$.
By the very definition, for any particular value of $Y$, we have
$
\E_X |\left<X,Y\right>|^p \leq M_p^p\,|Y|^p.
$
It remains to take the expectation with respect to $Y$.

In particular, $m_2 \leq M_2^2$, as can also be seen from (2.2).
The identities in (2.2) also show that, in the general non-isotropic case, $M_2^2$ 
may be larger than $m_2$.

Let us now turn to the functionals
$$
\sigma_{2p} = \sigma_{2p}(X) = 
\sqrt{n} \left(\E\, \Big|\frac{|X|^2}{n} - 1\Big|^p\right)^{1/p}, \qquad
p \geq 1,
$$
where it is natural to assume that $\E\,|X|^2 = n$.
Note that $\sigma_{2p}$ represents a non-decreasing function of $p$, which 
attains its minimum at $p=1$ with value
$$
\sigma_2 = \sigma_2(X) \, = \, \frac{1}{\sqrt{n}}\, \E\, \big|\,|X|^2 - n\big|.
$$
Another important value is $\sigma_4 = \frac{1}{n}\, \Var(|X|^2)$.
They may be related to the variance of the Euclidean norm.

\vskip5mm
{\bf Proposition 2.4.} {\sl If $\E\,|X|^2 = n$, then $\Var(|X|) \leq \sigma_4^2$.
In addition,
$$
\frac{1}{4}\,\sigma_2^2 \leq \Var(|X|) \leq \sigma_2\sqrt{n}.
$$
}

\vskip2mm
{\bf Proof.} Put $\xi = \frac{1}{\sqrt{n}}\,|X|$ and $a = \sqrt{\E \xi^2}$. Then,
since $\xi \geq 0$,
\bee
\Var(\xi^2) 
 & = &
\E\,(\xi^2 - a^2)^2 \\
 & = &
\E\,(\xi - a)^2(\xi + a)^2 \ \geq \ 
\E\,(\xi - a)^2 \cdot a^2 \ \geq \Var(\xi) \cdot a^2.
\ene
That is, $\E \xi^2 \, \Var(\xi) \leq \Var(\xi^2)$, which is exactly the first 
required relation.

Now, in terms of $\xi$, one may write
$$
\Var(|X|) = n\,\Var(\xi) = n\,(1 - (\E\xi)^2) = n\,(1 - \E \xi)\,(1 + \E \xi),
$$
while $\sigma_2 = \sqrt{n}\, \E\,|1 - \xi^2|$. By Cauchy's inequality,
$$
(\E\,|1 - \xi^2|)^2 \, \leq \, \E\,(1 - \xi)^2\,\E\,(1 + \xi)^2 \, = \,
4\, \E\,(1 - \xi)\,\E\,(1 + \xi),
$$
implying that $\sigma_2^2 \leq 4\,\Var(|X|)$.

The last inequality of the proposition may be rewritten as 
$1 - (\E \xi)^2 \leq \E\,|1 - \xi^2|$. If $(\Omega,\P)$ is the underlying 
probability space, define the probability measure 
$$
dQ = (1+\xi)\,d\P/\E\,(1 + \xi)
$$ 
and write $\E_Q$ for the expectation with respect to it. The required inequality
then takes the form $\E_Q\,|1 - \xi| \geq \E_Q (1-\xi)$, which is obvious.
\qed

\vskip5mm
The functionals $\sigma_{2p}^2$ and $m_p$ 
are useful in the problem of estimation of ``small" ball probabilities. 

\vskip5mm
{\bf Proposition 2.5.} {\sl Let $Y$ be an independent copy of a random vector 
$X$ in $\R^n$ such that $\E\, |X|^2 = n$. For all $p,q \geq 1$,
$$
\P\Big\{|X - Y|^2 \leq \frac{1}{4}\,n\Big\} \, \leq \, \frac{4^q}{n^{q/2}}\,m_q^q +
\frac{4^{2p}}{n^p}\,\sigma_{2p}^{2p}.
$$
In particular,
$$
\P\Big\{|X - Y|^2 \leq \frac{1}{4}\,n\Big\} \, \leq \, \frac{C}{n^p}
$$
with $C = 4^{2p}\,(m_{2p}^{2p} + \sigma_{2p}^{2p})$. 
}

\vskip5mm
{\bf Proof.} According to the definition,
$$
\sigma_{2p}^p = n^{-p/2}\, \E\, \big|\,|X|^2 - n\big|^p.
$$
Hence, for any $\lambda \in (0,1)$, by Chebyshev's inequality,
\bee
\P\big\{|X|^2 \leq \lambda n\big\}
 & = &
\P\Big\{\E\, |X|^2 - |X|^2 \geq (1-\lambda) \,\E\, |X|^2\Big\} \\
 & \leq &
\frac{\sigma_{2p}^p n^{p/2}}{(1-\lambda)^p\,(\E\, |X|^2)^p}
 \ = \ 
\frac{\sigma_{2p}^p}{(1-\lambda)^p\,n^{p/2}}.
\ene
In particular, choosing $\lambda = 3/4$, we get
$$
\P\Big\{|X|^2 + |Y|^2 \leq \frac{3}{4}\,n\Big\} \, \leq \,
\P\Big\{|X|^2 \leq \frac{3}{4}\,n\Big\}\,
\P\Big\{|Y|^2 \leq \frac{3}{4}\,n\Big\} \, \leq \,
\frac{4^{2p}\,\sigma_{2p}^{2p}}{n^p}.
$$ 

On the other hand, by Markov's inequality,
$$
\P\Big\{|\left<X,Y\right>| \geq \frac{1}{4}\,n\Big\}
 \leq \frac{4^q\, \E\, |\left<X,Y\right>|^q}{n^q} = \frac{4^q\,m_q^q}{n^{q/2}}.
$$
One may now write
$$
|X - Y|^2 = |X|^2 + |Y|^2 - 2 \left<X,Y\right>
$$
and split the event $|X - Y|^2 \leq \frac{1}{4}\,n$ into the case
$|\left<X,Y\right>| \geq \frac{1}{4}\,n$ and the case of the
opposite inequality. In view of the set inclusion
$$
\Big\{|X - Y|^2 \leq \frac{1}{4}\,n\Big\} \subset 
\Big\{|\left<X,Y\right>| \geq \frac{1}{4}\,n\Big\} \cup
\Big\{|X|^2 + |Y|^2 \leq \frac{3}{4}\,n\Big\},
$$
the proposition follows.
\qed


\vskip10mm
\section{{\bf Linear functionals on the sphere}}
\setcounter{equation}{0}

\vskip2mm
\noindent
The aim of this section is to quantify the  asymptotic normality of 
distributions of linear functionals with respect to the normalized Lebesgue 
measure $\mu_{n-1}$ on the unit sphere $S^{n-1} \subset \R^n$ $(n \geq 2)$. 
By the rotational invariance of this measure, all linear functionals 
$f(\theta) = \left<\theta,v\right>$ with $|v|=1$ have equal distributions, 
and it is sufficient to focus just on the first coordinate $\theta_1$ of 
the vector $\theta \in S^{n-1}$. As a random variable on the probability 
space $(S^{n-1},\mu_{n-1})$, it has density
$$
c_n\, \big(1 - x^2\big)_+^{\frac{n-3}{2}}, \qquad x \in \R,
$$
where 
$
c_n =  \frac{\Gamma(\frac{n}{2})}{\sqrt{\pi}\,\Gamma(\frac{n-1}{2})}
$
is a normalizing constant. 

Let us denote by $\varphi_n$ the density of the normalized first 
coordinate $Z_n = \sqrt{n}\, \theta_1$ under the measure $\mu_{n-1}$, i.e.,
$$
\varphi_n(x) =  c_n'\, \Big(1 - \frac{x^2}{n}\Big)_+^{\frac{n-3}{2}}, 
\qquad c_n' = \frac{c_n}{\sqrt{n}} = 
\frac{\Gamma\big(\frac{n}{2}\big)}{\sqrt{\pi n}\ \Gamma\big(\frac{n-1}{2}\big)}.
$$
Clearly, as $n \rightarrow \infty$,
$$
\varphi_n(x) \rightarrow \varphi(x) = \frac{1}{\sqrt{2\pi}}\,e^{-x^2/2},
\qquad
c_n'  \rightarrow \frac{1}{\sqrt{2\pi}},
$$
and one can show that $c_n'  < \frac{1}{\sqrt{2\pi}}$ for all $n \geq 2$.

We are interested in non-uniform deviation bounds of $\varphi_n(x)$ from 
$\varphi(x)$.

\vskip5mm
{\bf Proposition 3.1.} {\sl If $n \geq 3$, then for all $x \in \R$, with 
some universal constant $C$
\be
|\varphi_n(x) - \varphi(x)| \, \leq \, \frac{C}{n}\,e^{-x^2/8}.
\en
}

\vskip2mm
{\bf Proof.} Note that since the random variable $Z_3$ has uniform 
distribution on $[-\sqrt{3},\sqrt{3}]$, inequality (3.1) obviously holds 
for $n=3$. Hence, let $n \geq 4$.

First we consider the asymptotic behavior of the functions
$$
p_n(x) =  \Big(1 - \frac{x^2}{n}\Big)_+^{\frac{n-3}{2}}, \qquad x \in \R.
$$
Clearly, $p_n(x) \rightarrow e^{-x^2/2}$ for all $x$. Moreover, for 
$|x| < \sqrt{n}$, we have
$$
-\log p_n(x) \, = \, -\frac{n-3}{2}\,\log\Big(1 - \frac{x^2}{n}\Big)
 \, \geq \, \frac{n-3}{2}\,\frac{x^2}{n} \, \geq \, \frac{x^2}{8},
$$
so that there is a uniform bound
\be
p_n(x) \leq e^{-x^2/8}, \qquad x \in \R.
\en

To study the rate of convergence of $p_n(x)$, assume that 
$|x| \leq \frac{1}{2}\sqrt{n}$.
By Taylor's expansion, with some $0 \leq \ep \leq 1$
\bee
-\log p_n(x)
 & = &
\frac{n-3}{2}\,\Big[\,\frac{x^2}{n} +
\Big(\frac{x^2}{n}\Big)^2
\sum_{k=2}^\infty \frac{1}{k}\, \Big(\frac{x^2}{n}\Big)^{k-2}\,\Big] \\
 & = &
\frac{n-3}{2}\,\Big(\frac{x^2}{n} + \frac{x^4}{n^2}\,\ep\Big) 
 \, = \,
\frac{x^2}{2} + \frac{x^2}{2n}\, \Big(-3 + \frac{n-3}{n}\, x^2 \ep\Big),
\ene
that is,
$$
p_n(x) = e^{-x^2/2}\,e^{-\delta} \quad {\rm with} \ \ 
\delta = \frac{x^2}{2n}\,\Big(-3 + \frac{n-3}{n}\, x^2 \ep\Big).
$$
Since
$
\delta \geq -\frac{3x^2}{2n} \geq -\frac{3}{8n} \geq -\frac{3}{32},
$
we have
$$
|e^{-\delta} - 1| \leq |\delta|\,e^{3/32} \leq 1.1\,|\delta|.
$$
On the other hand,
$$
\delta \leq \frac{x^2}{2n}\,\Big(-3 + \frac{n-3}{n}\, x^2\Big)
\leq \frac{x^4}{2n}, 
$$
which together with the lower bound on $\delta$ yields
$$
1.1\,|\delta| \leq 1.1\,\Big(\frac{3x^2}{2n} + \frac{x^4}{2n}\Big) \leq
\frac{1}{n}\,(3x^2 + x^4).
$$

Thus,
$$
|p_n(x) - e^{-x^2/2}| \, \leq \, \frac{1}{n}\,(3x^2 + x^4)\, e^{-x^2/2}, \qquad 
|x| \leq \frac{1}{2}\sqrt{n}.
$$
Combining this inequality with (3.2), we also get a non-uniform
bound on the whole real line, namely
$$
|p_n(x) - e^{-x^2/2}| \, \leq \, \frac{C}{n}\,e^{-x^2/8}, \qquad x \in \R,
$$
where $C$ is an absolute constant. Let us integrate this inequality over $x$. Since
$$
\int_{-\infty}^{\infty} p_n(x)\,dx = \frac{1}{c_n'}, \qquad
\int_{-\infty}^{\infty} e^{-x^2/2}\,dx = \sqrt{2\pi},
$$
we get that $|\frac{1}{c_n'} - \sqrt{2\pi}| \leq \frac{C}{n}$
with some absolute constant $C$. Hence, we arrive at the conclusion (3.1) 
for the densities $\varphi_n$ for $n \ge 4$ as well.
\qed

\vskip5mm
In the sequel we denote by $J_n$ the characteristic function of the first 
coordinate $\theta_1$ of a random vector $\theta$ which is uniformly distributed on 
the unit sphere $S^{n-1}$. In a more explicit form, for any $t \in \R$,
\bee
J_n(t) 
 & = &
c_n \int_{-\infty}^{\infty} e^{itx}\,(1 - x^2)_+^{\frac{n-3}{2}}\,dx \\
 & = &
c_n' \int_{-\infty}^{\infty} e^{itx/\sqrt{n}}\,
\Big(1 - \frac{x^2}{n}\Big)_+^{\frac{n-3}{2}}\,dx.
\ene
Note that the equality
$$
\widetilde J_\nu(t) = 
\frac{1}{\sqrt{\pi}\,\Gamma(\nu + \frac{1}{2})}\,\Big(\frac{t}{2}\Big)^\nu
\int_{-1}^1 e^{itx}\,(1 - x^2)^{\nu - \frac{1}{2}}\,dx
$$
defines the classical Bessel function of the first kind with index $\nu$
([Ba], p.\,81). Therefore, 
$$
J_n(t) = \frac{1}{c_n}\sqrt{\pi}\,\Gamma(\nu + \frac{1}{2})\,
\Big(\frac{t}{2}\Big)^{-\nu}\,\widetilde J_\nu(t), \qquad \nu = \frac{n}{2} - 1.
$$
However, this relationship will not be used in the sequel.

Thus, the characteristic function of $Z_n = \theta_1 \sqrt{n}$ is given by
$$
\hat \varphi_n(t) = J_n\big(t\sqrt{n}\big) = 
\int_{-\infty}^\infty e^{itx}\varphi_n(x)\,dx,
$$
which is the Fourier transform of the probability density $\varphi_n$. 
One immediate consequence from Proposition 3.1 is the following:

\vskip5mm
{\bf Corollary 3.2.} {\sl For all $t \in \R$, we have
$$
\big|J_n\big(t\sqrt{n}\big) - e^{-t^2/2}\big| \leq \frac{C}{n},
$$
where $C$ is an absolute constant. 
}

\vskip5mm
For large $t$, this bound may be improved by virtue of the following upper bound.

\vskip5mm
{\bf Proposition 3.3.} {\sl For all $t \in \R$,
\be
\big|J_n(t\sqrt{n})\big|\, \leq \, 4.1\,e^{-t^2/2} + 4\, e^{-n/12}. 
\en
}

{\bf Proof.} One may assume that $n \geq 4$ (since $4\, e^{-n/12} > 1$
for $n = 2$ and $n=3$, while $|J_n| \leq 1$). 
For the approximation we shall use an approach based on contour integration 
in complex analyis.

The function $z \rightarrow (1 - z^2)^{\frac{n-3}{2}}$ is analytic in the 
whole complex plane when $n$ is odd and in the strip $z = x + iy$, $|x| < 1$, 
when $n$ is even. Therefore, integrating along the boundary of the rectangle 
${\mathfrak C} = [-1,1] \times [0,y]$ with $y > 0$ (slightly modifying the contour 
in a standard way near the points $-1$ and $1$), we have
$$
\int_{\mathfrak C} e^{itz}\, \big(1 - z^2\big)^{\frac{n-3}{2}}\,dz = 0.
$$
Then we obtain a natural decomposition
$
J_n(t\sqrt{n}) = c_n\, \big(I_1(t) + I_2(t) + I_3(t)\big) 
$
for $t>0$, where
$$
I_1(t) \ = \ e^{-ty\sqrt{n}} \, \big(1+y^2\big)^{\frac{n-3}2}
\int_{-1}^1 e^{itx\sqrt{n}}\,
\Big(\frac{1-(x+iy)^2}{1+y^2}\Big)^{\frac{n-3}2}\,dx,
$$
$$
I_2(t) \ = \ - e^{it\sqrt n} \int_0^y e^{-ts\sqrt n} \, 
\big(1-(1 + is)^2\big)^{\frac{n-3}2}\,ds,
$$
$$
I_3(t) \ = \ e^{-it\sqrt n} \int_0^y e^{-ts\sqrt n} \, 
\big(1-(1 - is)^2\big)^{\frac{n-3}2}\,ds.
$$
For $0 \le s \leq y\le \alpha$,
$$
\big|\,1-(1 + is)^2\big| \, = \,s\sqrt{s^2+4} \, \leq \,
\alpha\, \sqrt{\alpha^2 + 4} \equiv \beta.
$$
Choosing $\alpha = \frac{1}{\sqrt{6}}$, we have $\beta = \frac{5}{6}$.
Hence, for all $t>0$,
\be
|I_2(t)| \, \le \,\beta^{\frac{n-3}{2}} \, 
\int_0^y e^{-ts\sqrt{n}}\,ds \, \le \, \frac{1}{t\sqrt{n}}\, \beta^{\frac{n-3}{2}}.
\en

The same estimates hold for $I_3(t)$. 

In order to estimate $I_1(t)$, we use an elementary identity
$$
\big|1-(x+iy)^2\big|^2 = (1-2x^2 )(1+y^2)^2 + x^2\,(x^2 + 6y^2 + 2y^4) \qquad 
(x,y \in \R),
$$
which for the region $|x| \leq 1$ yields
$$
\bigg|\frac{1-(x+iy)^2}{1+y^2}\bigg|^2 \, \le \, 
1 - 2x^2 + x^2\,v(y^2), \qquad v(z) = \frac{1 + 6z + 2z^2}{(1 + z)^2}.
$$
Since $v'(z) = \frac{4-2z}{(1+z)^3} > 0$, this function increases in 
$0 \leq z \leq 2$, and since $z = y^2 \leq \frac{1}{6}$, we have 
$v(y^2) \leq v(1/6) = \frac{74}{49}$. Hence
$$
\bigg|\frac{1-(x+iy)^2}{1+y^2}\bigg|^2 \, \le \, 1 - \frac{24}{49}\,x^2 \, \leq \, 
e^{-\frac{24}{49}\,x^2}.
$$
Using this estimate together with $n - 3 \geq \frac{1}{4}\,n$, we have
\bee
\int_{-1}^1 \Big(\frac{|1-(x + iy)^2|}{1+y^2}\Big)^{\frac{n-3}2}\,dx 
 & \le &
\int_1^1 e^{-\frac{12}{49}\,\frac{n-3}{2}\,x^2}\,dx \\
 & = & 
\sqrt{\frac{49}{12}} \, \sqrt{\frac{2\pi}{n-3}} \, \leq \, 
\frac{14}{\sqrt{12}}\,\sqrt{\frac{2\pi}{n}}
 \, \leq \, 4.1\,\sqrt{\frac{2\pi}{n}}.
\ene
This upper bound allows us to conclude that
\bee
|I_1(t)|
 & \le & 
4.1\,\sqrt{\frac{2\pi}{n}} \, e^{-ty\sqrt{n}} \, \big(1+y^2\big)^{\frac{n-3}{2}} \\
 & \le &
4.1\,\sqrt{\frac{2\pi}{n}} \, \exp\Big\{-ty\sqrt{n} + \frac{n-3}{2}\,y^2\Big\}.
\ene
Choosing here $y = \frac{t}{\sqrt{n}}$, the expression in the exponent will 
be smaller than $t^2/2$, hence
\be
|I_1(t)| \le 4.1\, e^{-t^2/2}, \qquad 0 \leq t \leq \sqrt{n/6}.
\en
In the case $t > \sqrt{n/6}$, we choose $y = \frac{1}{\sqrt{6}}$ and then
$-ty\sqrt{n} + \frac{n-3}{2}\,y^2 < - \frac{n}{12} - \frac{1}{4}$, so that
\be
|I_1(t)| \le 4.1\,\sqrt{\frac{2\pi}{n}}\, e^{- \frac{n}{12} - \frac{1}{4}}, 
\qquad t \geq \sqrt{n/6}.
\en

Let us collect these estimates. For $2 \leq t \leq \sqrt{n/6}$, we combine 
(3.5) with (3.4) and a similar bound for $I_3(t)$, and use 
$c_n < \frac{1}{\sqrt{2\pi}} \sqrt{n}$ with
$\beta^{\frac{n-3}{2}} < 0.77\,e^{-n/12}$. This leads to
\bee
c_n\, \big(|I_1(t)| + |I_2(t)| + |I_3(t)|\big) 
 & \leq & 
\frac{c_n}{\sqrt{n}}\, \beta^{\frac{n-3}{2}} + 
4.1\,c_n \sqrt{\frac{2\pi}{n}}\, e^{-t^2/2} \\
 & \leq &
0.31\,e^{-n/12} + 4.1\,c_n \, e^{-t^2/2}.
\ene
Similarly, in case $t > \sqrt{n/6}$, we use (3.6) leading to
$$
c_n\, \big(|I_1(t)| + |I_2(t)| + |I_3(t)|\big) \leq 
0.31\,e^{-n/12} + 4.1\,e^{-1/4}\,e^{-n/12} < 4 e^{-n/12}.
$$

Finally, if $0 \leq t \leq 2$, then $|J_n(t\sqrt{n})| \leq 1 \leq 4\, e^{-t^2/2}$.
\qed


\vskip10mm
\section{{\bf Typical Distributions and Mixtures of Gaussian Measures}}
\setcounter{equation}{0}

\vskip2mm
\noindent
The asymptotic normality of the typical distributions $F$ in Sudakov's theorem, 
defined in (1.1), may be described in the next assertion proved in [B-C-G1].

\vskip5mm
{\bf Proposition 4.1.} {\sl Given a random vector $X$ in $\R^n$,
suppose that $\E\, |X|^2 = n$. With some absolute constant $c>0$ we have
\be
\int_{-\infty}^\infty (1 + x^2)\,|F(dx) - \Phi(dx)| \, \leq \, 
c\,\Big(\frac{1}{n} + \Var(r)\Big),
\en
where $r = \frac{1}{\sqrt{n}} |X|$.
}

\vskip5mm
Here the positive measure $|F - \Phi|$ denotes the variation in the sense of 
measure theory, and the left integral represents the weighted total variation 
of $F - \Phi$. In particular, we have a similar bound for the usual total 
variation distance between $F$ and $\Phi$, as well as for the Kolmogorov distance 
$\rho(F,\Phi)$. Applying Proposition 2.4, the latter
may be related to the variance-type functionals $\sigma_{2p}$ (cf. also [M-M]).

\vskip5mm
{\bf Corollary 4.2.} {\sl In particular (under the same conditions),
$$
\rho(F,\Phi) \, \leq \, c\,\frac{1 + \sigma_4^2}{n}, \qquad
\rho(F,\Phi)  \, \leq \, c\,\frac{1 + \sigma_2}{\sqrt{n}}.
$$
}

The proof of Proposition 4.1 is based on the following observation about 
general mixtures of centered Gaussian measure on the real line. Given a random 
variable $r \geq 0$, let us denote by $\Phi_r$ the distribution function 
of the random variable $r Z$, where $Z \sim N(0,1)$ is independent of $r$. That is,
$$
\Phi_r(x) = \P\{r Z \leq x\} = \E\,\Phi(x/r), \qquad x \in \R.
$$
As  shown in [B-C-G1], if $\E r^2 = 1$, then with some absolute constant 
$c$ we have
\be
\int_{-\infty}^\infty (1 + x^2)\, |\Phi_r - \Phi|(dx) \, \leq \, c\,\Var(r).
\en

To explain the transition from (4.2) to (4.1), assume that $n \geq 3$.
Let $\Phi_n$ and $\varphi_n$ denote respectively the distribution function and 
the density of $Z_n = \theta_1 \sqrt{n}$, where $\theta_1$ is the first coordinate 
of a random point $\theta$ uniformly distributed in $S^{n-1}$. If 
$r^2 = \frac{1}{n}\,|X|^2$ is independent of $Z_n$ ($r \ge 0$), then, by 
the definition of the typical distribution,
$$
F(x)  =  \P\{r Z_n \leq x\} \, = \, \E\,\Phi_n(x/r), \qquad x \in \R,
$$
so that
\be
\int_{-\infty}^\infty (1 + x^2)\,|F(dx) - \Phi_r(dx)| \, = \,
\int_{-\infty}^\infty (1 + x^2)\,|\E\,\Phi_n(dx/r) - \E\,\Phi(dx/r)|.
\en
But, for any fixed value of $r$,
$$
\int_{-\infty}^\infty (1 + x^2)\,|\Phi_n(dx/r) - \Phi(dx/r)| \, = \,
\int_{-\infty}^\infty (1 + r^2 x^2)\,|\Phi_n(dx) - \Phi(dx)|,
$$
hence, by (4.2), taking the expectation with respect to $r$ and using 
Jensen's inequality, we get
\bee
\int_{-\infty}^\infty (1 + x^2)\,|F(dx) - \Phi_r(dx)|
 & \leq &
\E \int_{-\infty}^\infty (1 + x^2)\,|\Phi_n(dx/r) - \Phi(dx/r)| \\ 
 & = &
\E \int_{-\infty}^\infty (1 + r^2 x^2)\,|\Phi_n(dx) - \Phi(dx)| \\ 
 & = &
\int_{-\infty}^\infty (1 + x^2)\,|\Phi_n(dx) - \Phi(dx)|.
\ene
It remains to apply (3.1), which yields
$$
\int_{-\infty}^\infty (1 + x^2)\,|\Phi_n(dx) - \Phi(dx)| \, = \,
\int_{-\infty}^\infty (1 + x^2)\,|\varphi_n(x) - \varphi(x)|\,dx
 \, \leq \, \frac{C}{n}
$$
with some universal constant $C$.
\qed


\vskip10mm
\section{{\bf Characteristic Functions of Weighted Sums}}
\setcounter{equation}{0}

\vskip2mm
\noindent
As before, let $X = (X_1,\dots,X_n)$ denote a random vector in $\R^n$, 
$n \geq 2$. The concentration problems for distributions of weighted sums 
$S_\theta = \left<X,\theta\right>$ may be studied by means of their 
characteristic functions
\be
f_\theta(t) = \E\,e^{it \left<X,\theta\right>},\qquad t \in \R.
\en
In particular, we intend to quantify the concentration of $f_\theta$ around the 
characteristic function $f$ of the typical distribution $F$ on average over 
the directions $\theta$ in terms of correlation-type functionals. 
Note that the characteristic function of $F$ is given by
$$
f(t) = \E_\theta f_\theta(t) = \E_\theta \,
\E\,e^{it \left<X,\theta\right>} = \E\, J_n(t|X|), \qquad t \in \R,
$$
where $J_n$ is the characteristic function of the first coordinate 
$\theta_1$ under the uniform measure $\mu_{n-1}$ on the unit sphere $S^{n-1}$.

First let us describe the decay of $t \rightarrow |f_\theta(t)|$ at infinity 
on average with respect to $\theta$. Starting from (5.1), write
$$
\E_\theta\, |f_\theta(t)|^2 = \E_\theta\, e^{it \left<X-Y,\theta\right>} = 
\E\, J_n(t|X-Y|),
$$
where $Y$ is an independent copy of $X$. To proceed, let us rewrite the Gaussian-type 
bound (3.3) of Proposition 3.3 as
\be
|J_n(t)|\, \leq \,4.1\,e^{-t^2/2n} + 4\,e^{-n/12}
\en
which gives
$$
\E_\theta\, |f_\theta(t)|^2\, \leq \, 4.1\,\E\,e^{-t^2 |X-Y|^2/2n} + 4\,e^{-n/12}.
$$
Splitting the latter expectation into the event $A = \{|X-Y|^2 \leq \lambda n\}$ 
and its complement, we get the following general bound.

\vskip5mm
{\bf Lemma 5.1.} {\sl The characteristic functions $f_\theta$ satisfy, for all 
$t \in \R$ and $\lambda > 0$,
$$
\frac{1}{2.1}\,\E_\theta\, |f_\theta(t)| \leq e^{-\lambda t^2/4} + e^{-n/24} + 
\sqrt{\P\{|X-Y|^2 \leq \lambda n\}}\,,
$$
where $Y$ is an independent copy of $X$.
}

\vskip5mm
In case $\E\, |X|^2 = n$, the right-hand side of these bounds can be further 
quantified by using the moment and variance-type functionals, which we have 
discussed before, namely
$$
m_p = \frac{1}{\sqrt{n}}\,(\E\, |\left<X,Y\right>|^p)^{1/p}, \qquad
\sigma_{2p} = \sqrt{n} \left(\E\, \Big|\frac{|X|^2}{n} - 1\Big|^p\right)^{1/p}.
$$
Note that both $m_p$ and $\sigma_p$ are non-decreasing functions in $p \geq 1$. 
In order to estimate the probability of the event $B$, we shall use Proposition 2.5,
which gives
$$
\P(A) \, \leq \, \frac{C}{n^p}
$$
with a constant $C = 4^{2p}\,(m_{2p}^{2p} + \sigma_{2p}^{2p})$.
Hence, from Lemma 5.1 and using $m_{2p} \geq m_2 \geq 1$ (cf. Proposition 2.1), 
we deduce:

\vskip5mm
{\bf Lemma 5.2.} {\sl Suppose that $\E\, |X|^2 = n$. If the moment $m_{2p}$ is 
finite for $p \geq 1$, then with some constant $c_p>0$ depending on $p$,
$$
c_p\,\E_\theta\, |f_\theta(t)| \, \leq \,
\frac{m_{2p}^p + \sigma_{2p}^p}{n^{p/2}} + e^{-t^2/16}.
$$
}

\vskip2mm
By the triangle inequality, $|f(t)| \leq \E_\theta\, |f_\theta(t)|$. Hence, 
the characteristic function of the typical distribution shares the same bounds. 
In fact, here the parameter $m_{2p}$ is not needed. Indeed, as was shown in 
the proof of Proposition 2.5 with $\lambda = \frac{1}{2}$, we have
$$
\P\Big\{|X|^2 \leq \frac{1}{2}\, n\Big\} \,\leq\, 2^p \, 
\frac{\sigma_{2p}^p}{n^{p/2}}.
$$ 
Hence, by (5.2),
\bee
|f(t)| 
 & \leq &
\E\, |J_n(t|X|)|\,1_{\{|X| \leq \sqrt{n/2}\}} +
\E\, |J_n(t|X|)|\,1_{\{|X| > \sqrt{n/2}\}} \\
 & \leq &
2^p \, \frac{\sigma_{2p}^p}{n^{p/2}} + C\, \big(e^{-t^2/4} + e^{-n/12}\big).
\ene
Thus, we get:

\vskip5mm
{\bf Lemma 5.3.} {\sl Suppose that $\E\, |X|^2 = n$. Then with some constant 
$c_p>0$ depending on $p \geq 1$, for all $t \in \R$,
$$
c_p\,|f(t)| \leq \frac{1 + \sigma_{2p}^p}{n^{p/2}} + e^{-t^2/4},
$$
and therefore, for all $T>0$,
$$
\frac{c_p}{T} \int_0^T |f(t)|\,dt \, \leq \, 
\frac{1 + \sigma_{2p}^p}{n^{p/2}} + \frac{1}{T}.
$$
}

We first study the concentration properties of $f_\theta(t)$ as functions of 
$\theta$ on the sphere with fixed $t \in \R$ (rather than directly for 
the distributions $F_\theta$). This can be done in terms of the moment functionals
$$
M_p = M_p(X) \, = 
\sup_{\theta \in S^{n-1}} \Big(\E\, |\left<X,\theta\right>|^p\Big)^{1/p}.
$$

Our basic tool is a well-known spherical Poincar\'e inequality
\be
\int_{S^{n-1}} |u(\theta) - a|^2\,d\mu_{n-1}(\theta) \leq 
\frac{1}{n-1}\,\int_{S^{n-1}} |\nabla u(\theta)|^2\,d\mu_{n-1}(\theta).
\en
It holds true for any complex-valued function $u$ which is defined and smooth 
in a neighborhood of the sphere, and has gradient $\nabla u$ and the mean 
$a = \int u\,d\mu_{n-1}$ (cf. [L]).

According to (5.1), the function $\theta \rightarrow f_\theta(t)$ is smooth on 
the whole space $\R^n$ and has partial derivatives
$$
\frac{\partial_j f_\theta(t)}{\partial \theta_j} = 
it\, \E X_j\, e^{it \left<X,\theta\right>}
$$
or in the vector form
$$
\left<\nabla f_\theta(t),v\right> = 
it\, \E \left<X,v\right> e^{it \left<X,\theta\right>}, \qquad v \in \R^n.
$$
Hence
$$
|\left<\nabla f_\theta(t),v\right>| \leq |t|\, \E\, |\left<X,v\right>|.
$$
Taking the sup over all $v \in S^{n-1}$, we obtain a uniform bound on the
modulus of the gradient, namely $|\nabla f_\theta(t)| \leq M_1 |t|$.

A similar bound holds as well in average. To this aim, let us square 
the vector representation and write
$$
\left<\nabla f_\theta(t),v\right>^2 = 
t^2\, \E \left<X,v\right>\left<Y,v\right> e^{it \left<X-Y,\theta\right>}, 
$$
where $Y$ is an independent copy of $X$. Integrating over $v$ with respect
to $\mu_{n-1}$, we get the representation
$$
|\nabla f_\theta(t)|^2 \, = \, 
t^2\, \E \left<X,Y\right> e^{it \left<X-Y,\theta\right>},
$$
so that
$$
\E_\theta\,|\nabla f_\theta(t)|^2 \, = \, t^2\, \E \left<X,Y\right> J_n(t (X-Y)).
$$
(where $\E_\theta$ refers to integration over $\mu_{n-1}$).
Applying (5.3), one can summarize.

\vskip5mm
{\bf Lemma 5.4.} {\sl Given a random vector $X$ in $\R^n$ with finite moment 
$M_1$, for all $t \in \R$,
$$
\E_\theta\,|f_\theta(t) - f(t)|^2  \leq \frac{t^2}{n-1}\,M_1^2.
$$
In addition,
$$
\E_\theta\,|f_\theta(t) - f(t)|^2  \leq 
\frac{t^2}{n-1}\,\E \left<X,Y\right> J_n(t (X-Y)),
$$
where $Y$ is an independent copy of $X$.
}


\vskip10mm
\section{{\bf Berry-Esseen Bounds. Theorem 1.2 and its Generalization}}
\setcounter{equation}{0}

\vskip2mm
\noindent
Fourier Analysis provides a well-established tool to prove Berry-Esseen-type 
bounds for the Kolmogorov distance
$$
\rho(F_\theta,F) = \sup_x\, |F_\theta(x) - F(x)|.
$$
To study the average behavior of this distance with respect to $\theta$ using 
the uniform measure $\mu_{n-1}$ on the unit sphere, as a preliminary step, 
let us first introduce two auxiliary bounds.

\vskip5mm
{\bf Lemma 6.1.} {\sl Let $X$ be a random vector in $\R^n$. With some absolute 
constant $c>0$, for all $T \geq T_0 > 0$,
\begin{eqnarray}
c\,\E_\theta\, \rho(F_\theta,F) 
 &  \leq & 
\int_0^{T_0} \frac{\E_\theta\,|f_\theta(t) - f(t)|}{t}\,dt \nonumber \\ 
 & & + 
\int_{T_0}^T \frac{\E_\theta\,|f_\theta(t)|}{t}\,dt + \frac{1}{T} 
\int_0^T |f(t)|\,dt.
\end{eqnarray}
}

As before, here  $F_\theta$ denote distribution functions of the weighted sums
$S_\theta = \left<X,\theta\right>$ with their characteristic functions
$$
f_\theta(t) = \E\,e^{it\left<X,\theta\right>} = 
\int_{-\infty}^\infty e^{itx}\, dF_\theta(x), \qquad t \in \R, \ \theta \in S^{n-1},
$$
and $F = \E_\theta F$ is the typical distribution function with characteristic 
function
$$
f(t) = \E_\theta f_\theta(t)= \int_{-\infty}^\infty e^{itx}\, dF(x).
$$

For an estimation of the Kolmogorov distance, the following general Berry-Esseen 
bound will be convenient:
\be
c\,\rho(U,V) \leq \int_0^T \frac{|u(t) - v(t)|}{t}\,dt + \frac{1}{T} 
\int_0^T |v(t)|\,dt \qquad (T>0).
\en
Here $U$ and $V$ may be arbitrary distribution functions on the line 
with characteristic functions $u$ and $v$, respectively, and $c>0$ is 
an absolute constant (cf. e.g. [P1-2], [B3]).

In our situation, we take $U = F_\theta$ and $V = F$.
In order to estimate the first integral in (6.2), we shall split the integration 
into the two intervals, $[0,T_0]$ (the interval of moderate values of $t$), where 
it is easier to control the closeness of the two characteristic functions, and 
the long interval $[T_0,T]$, where both characteristic functions can be shown 
to be sufficiently small. Note that, by the triangle inequality, we have 
$|f(t)| \leq \E_\theta\, | f_\theta(t)|$, which implies 
$\E_\theta\,|f_\theta(t) - f(t)| \leq 2\,\E_\theta\,|f_\theta(t)|$.
Using this on the long interval, we arrive at the more specific variant of (6.2), 
namely (6.1).

The estimation of the integrals in (6.1) will be done in terms of
the functionals $m_p = m_p(X)$, $M_p = M_p(X)$ and $\sigma_{2p} = \sigma_{2p}(X)$.

\vskip5mm
{\bf Lemma 6.2.} {\sl Suppose that $X$ has \tc{a} finite moment of order $2p$ 
$(p \geq 1)$, and $\E\,|X|^2 = n$. Then with some 
constant $c_p$ depending on $p$ only, for all $T \geq T_0 > 0$,
\bee
c_p\,\E_\theta\, \rho(F_\theta,F)
 & \leq & 
\int_0^{T_0} \E_\theta\,|f_\theta(t) - f(t)|\,\frac{dt}{t} \\
 & & + \ 
\frac{m_{2p}^p + \sigma_{2p}^p}{n^{p/2}}\,\bigg(1 + \log\frac{T}{T_0}\bigg)
+ \frac{1}{T} + e^{-T_0^2/16}.
\ene
}

{\bf Proof.} By the second inequality of Lemma 5.3 
(on this step we use the assumption $\E\, |X|^2 = n$), we have
$$
\frac{c_p}{T} \int_0^T |f(t)|\,dt \, \leq \, \frac{1 + \sigma_{2p}^p}{n^{p/2}} + 
\frac{1}{T},
$$
while by Lemma 5.2 yields the bound
$$
c_p \int_{T_0}^T \frac{\E_\theta\, |f_\theta(t)|}{t}\,dt \, \leq \, 
\frac{m_{2p}^p + \sigma_{2p}^p}{n^{p/2}}\,\log\frac{T}{T_0} + e^{-T_0^2/16}.
$$
This allows us to estimate the pre-last and last integrals in (6.1).
\qed

\vskip5mm
We are now prepared to establish Theorem 1.2, in fact -- in somewhat more 
general form which requires the first moment, only. Recall that
$$
\sigma_2 = \frac{1}{\sqrt{n}}\, \E\, \big|\,|X|^2 - n\big|.
$$

\vskip5mm
{\bf Theorem 6.3.} {\sl If the random vector $X$ in $\R^n$ has finite first 
moment $M_1$, then
\be
c\,\E_\theta\, \rho(F_\theta,F) \, \leq \, M_1\sqrt{\frac{\log n}{n}} + 
\sqrt{\P\big\{|X-Y|^2 \leq n/4\big\}}\, \log n + \frac{1}{n},
\en
where $c>0$ is an absolute constant, and $Y$ is an independent copy of $X$.
As a consequence, if $X$ has finite 2-nd moment $M_2$ and $\E\, |X|^2 = n$, 
then
\be
c\,\E_\theta\, \rho(F_\theta,F) \, \leq \, 
(M_1 + m_2 + \sigma_2)\,\frac{\log n}{\sqrt{n}}.
\en
A similar bound also holds for the normal distribution function $\Phi$ 
in place of $F$.
}

\vskip5mm
The coefficient in (6.4) may be simplified by using $m_2 \leq M_2^2$ and 
$M_1 \leq M_2$. Since necessarily $M_2 \geq 1$, (6.4) implies the inequality 
(1.3) of Theorem 1.2.

\vskip5mm
{\bf Proof.} We apply Lemma 6.1 with $T_0 = 5\sqrt{\log n}$ and $T = 5n$.
The first integral in (6.1) can be bounded by virtue of the spherical 
Poincar\'e-type inequality, i.e., using the first bound of Lemma 5.4. It gives
$$
\E_\theta\,|f_\theta(t) - f(t)| \leq \frac{M_1 t}{\sqrt{n-1}} \qquad (t \geq 0)
$$
and hence
$$
\int_0^{T_0} \frac{\E_\theta\,|f_\theta(t) - f(t)|}{t}\,dt \, \leq \,
\frac{5 M_1}{\sqrt{n-1}} \sqrt{\log n}.
$$ 

Next, we apply Lemma 5.1 with $\lambda = 1/4$ which gives
\bee
\frac{1}{T} \int_0^T |f(t)|\,dt
 & \leq & 
\frac{1}{T} \int_0^T \E_\theta\, |f_\theta(t)|\,dt \\
 & \leq & 
\frac{2.1}{T} \int_0^T \Big(e^{-t^2/16} + e^{-n/24} + 
\sqrt{\P\{|X-Y|^2 \leq n/4\}}\,\Big)\,dt \\
 & \leq & 
\frac{c}{T} + 2.1 \sqrt{\P\{|X-Y|^2 \leq n/4\}}
\ene
with some absolute constant $c>0$. Similarly, 
$$
c\,\int_{T_0}^T \frac{\E_\theta\,|f_\theta(t)|}{t}\,dt \, \leq \, 
\Big(e^{-n/24} + \sqrt{\P\{|X-Y|^2 \leq n/4\}}\,\Big) \log\frac{T}{T_0}
+ e^{-T_0^2/16}.
$$
These bounds prove the first assertion of the theorem.

For the second assertion,  it remains to recall that, by Proposition 2.5, 
$$
\P\{|X-Y|^2 \leq n/4\} \leq 16\,\frac{m_2^2 + \sigma_2^2}{n},
$$
so that from (6.3) we get
\be
c\,\E_\theta\, \rho(F_\theta,F) \, \leq \, M_1\sqrt{\frac{\log n}{n}} + 
4\,\frac{m_2 + \sigma_2}{\sqrt{n}}\, \log n + \frac{1}{n}.
\en
Here, the last term $1/n$ is dominated by $m_2/n$.
This leads to the bound (6.4), in which $F$ may be replaced with 
the standard normal distribution function $\Phi$ due to the estimate
$\rho(F,\Phi) \leq \frac{C}{\sqrt{n}}\,\big(1 + \sigma_2\big)$, cf. Corollary 4.2.
\qed

\vskip5mm
{\bf Remark.} 
Working with the L\'evy distance $L$, which in general is weaker then the 
Kolmogorov distance $\rho$, one can get guaranteed rates with respect to $n$ for 
$\E_\theta\, L(F_\theta,F)$ in terms of $M_1$ or $M_2$. In particular, if 
$X$ isotropic, it is known that
$$
\mu_{n-1}\{L(F_\theta,F) \geq \delta\} \leq 4n^{3/8}\,e^{-n\delta^4/8}, \qquad
\delta > 0.
$$
This deviation bound yields
$$
\E_\theta\, L(F_\theta,F) \leq C\,\Big(\frac{\log n}{n}\Big)^{1/4}
$$
with some absolute constant $C$ ([B1]). See also [B2] for similar results about 
the Kantorovich distance.


\vskip10mm
\section{{\bf Proof of Theorem 1.1}}
\setcounter{equation}{0}

\vskip2mm
\noindent
In order to get rid of the logarithmic term in the bounds of Theorems 1.2/6.3, 
one may involve the 3-rd moment assumptions in terms of the moment 
and variance-type functionals $m_p$ and $\sigma_p$ of index $p=3$. They
are defined by
$$
m_3 = m_3(X) = \frac{1}{\sqrt{n}}\,\big(\E\, |\left<X,Y\right>|^3\big)^{1/3},
$$
where $Y$ is an independent copy of $X$, and
$$
\sigma_3 = \sigma_3(X) = 
\sqrt{n} 
\left(\E\, \Big|\frac{|X|^2}{n} - 1\Big|^{\frac{3}{2}}\right)^{\frac{2}{3}} = 
\frac{1}{\sqrt{n}} 
\left(\E\, \big|\,|X|^2 - n\big|^{\frac{3}{2}}\right)^{\frac{2}{3}}.
$$
Let us recall that $m_3 \leq M_3^2$. Hence, Theorem 1.1 will follow from 
the following, slightly sharpened assertion.

\vskip5mm
{\bf Theorem 7.1.} {\sl Let $X$ be a random vector in $\R^n$ with finite 3-rd 
moment, and such that $\E\,|X|^2 = n$. Then with some absolute constant $c$
\be
\E_\theta\, \rho(F_\theta,\Phi) \, \leq \, 
c\,(m_3^{3/2} + \sigma_3^{3/2})\,\frac{1}{\sqrt{n}}.
\en
}

\vskip2mm
{\bf Proof.} We now apply Lemma 6.2, choosing there $p = 3/2$, $T = 4n$ and 
$T_0 = 4\sqrt{\log n}$. Since necessarily $m_3 \geq 1$, the last term 
$e^{-T_0^2/16}$ is negligible, and we get the bound
$$
c\,\E_\theta\, \rho(F_\theta,F) \, \leq \, 
\int_0^{T_0} \frac{\E_\theta\,|f_\theta(t) - f(t)|}{t}\,dt + 
(m_3^{3/2} + \sigma_3^{3/2})\,\frac{\log n}{n^{3/4}}
$$
with some absolute constant $c>0$.
To analyze the last integral over the interval $[0,T_0]$, we apply Lemma 5.4, 
which gives
$$
\E_\theta\,|f_\theta(t) - f(t)| \, \leq \, \frac{t}{\sqrt{n-1}}\, 
\sqrt{\E \left<X,Y\right> J_n(t (X-Y))}, \qquad t \geq 0, 
$$
and hence
\begin{eqnarray}
c\,\E_\theta\, \rho(F_\theta,F) 
 & \leq & 
\frac{1}{\sqrt{n}} \int_0^{T_0} \sqrt{\E \left<X,Y\right> J_n(t (X-Y))}\,dt 
\nonumber \\ 
 & & + \  
(m_3^{3/2} + \sigma_3^{3/2})\,\frac{\log n}{n^{3/4}}.
\end{eqnarray}

Next, let us apply the bound of Corollary 3.2,
$|J_n\big(t\sqrt{n}\big) - e^{-t^2/2}| \leq \frac{C}{n}$, which allows one to replace
the $J_n$-term with $e^{-t^2 |X-Y|^2/2n}$ at the expense of an error of order
$$
\frac{1}{n}\,T_0 \sqrt{\E\,|\left<X,Y\right>|} \leq 
\frac{\sqrt{m_2}}{n^{3/4}}\,\,T_0 \leq 
\frac{m_3^{3/2}}{n^{3/4}}\,\,T_0,
$$
where we used the inequality $m_3 \geq m_2 \geq 1$.
As a result, the bound (7.2) may be simplified to
\be
c\,\E_\theta\, \rho(F_\theta,F) \, \leq \, 
\frac{1}{\sqrt{n}} \int_0^{T_0} \sqrt{I(t)}\,dt + (m_3^{3/2} + \sigma_3^{3/2})\,
\frac{\log n}{n^{3/4}}
\en
with 
$$
I(t) = \E \left<X,Y\right> e^{-t^2 |X-Y|^2/2n}.
$$
Note that $I(t) \geq 0$ which follows from
$I(t) = \big|\,\E\,e^{it\left<X,Z\right>/\sqrt{n}}\,|^2$, where the random vector
$Z$ is independent of $X$ and has a standard normal distribution on $\R^n$. 

Now, focusing on $I(t)$, consider the events 
$$
A = \Big\{|X-Y|^2 \leq \frac{1}{4}\,n\Big\}, \qquad
B = \Big\{|X-Y|^2 > \frac{1}{4}\,n\Big\}.
$$
We split the expectation in the definition of $I(t)$ into the sets
$A$ and $B$, so that $I(t) = I_1(t) + I_2(t)$, where
$$
I_1(t) = \E \left<X,Y\right> e^{-t^2 |X-Y|^2/2n}\,1_A, \qquad
I_2(t) = \E \left<X,Y\right> e^{-t^2 |X-Y|^2/2n}\,1_B.
$$
As we know (cf. Proposition 2.5),
$$
\P(A) \, \leq \, 64\,\frac{m_3^3 + \sigma_3^3}{n^{3/2}}.
$$
Hence, applying H\"older's inequality, we have
\bee
|I_1(t)| 
 & \leq & 
\big(\E\,|\left<X,Y\right>|^3)^{1/3}\ (\P(A))^{2/3}  \\
 & = &
m_3 \sqrt{n}\cdot 16\, \frac{m_3^2 + \sigma_3^2}{n}
\, \leq \, \frac{32}{\sqrt{n}}\,(m_3^3 + \sigma_3^3),
\ene
where we used that $m_3 \geq 1$.

Now, we represent the second expectation as
\bee
I_2(t) 
 & = &
e^{-t^2}\,
\E \left<X,Y\right> e^{-t^2 \big(\frac{|X-Y|^2}{2n} - 1\big)}\,1_B \\
 & = &
e^{-t^2}\,
\E \left<X,Y\right> 
\Big(e^{-t^2 \big(\frac{|X-Y|^2}{2n} - 1\big)} - 1\Big)\,1_B
- e^{-t^2}\,\E \left<X,Y\right> \,1_A.
\ene
Here the last expectation has been already bounded by 
$\frac{32}{\sqrt{n}}\, (m_3^3 + \sigma_3^3)$.
To estimate the first one, we use an elementary inequality
$$
|e^{-x} - 1| \leq |x|\,e^{x_0} \qquad (x_0 \geq 0, \ \ x \geq -x_0).
$$
Since on the set $B$, there is a uniform bound
$t^2 (\frac{|X-Y|^2}{2n} - 1) \geq -\frac{7}{8}\,t^2$, we conclude by virtue of
H\"older's inequality that
\bee
\E\, |\left<X,Y\right>|\, 
\big|\,e^{-t^2 \big(\frac{|X-Y|^2}{2n} - 1\big)} - 1\big|\,1_B 
 & \leq &
t^2 e^{7t^2/8}\,
\E\, |\left<X,Y\right>|\, \Big|\frac{|X-Y|^2}{2n} - 1\Big| \\ 
 & \hskip-20mm \leq & \hskip-10mm
t^2 e^{7t^2/8}\, \Big(\E\, |\left<X,Y\right>|^3\Big)^{\frac{1}{3}}
\bigg(\E\,\Big|\frac{|X-Y|^2}{2n} - 1\Big|^{\frac{3}{2}}\bigg)^{\frac{2}{3}}.
\ene
The first expectation on the right-hand side is
$\E\, |\left<X,Y\right>|^3 = (m_3\,\sqrt{n})^3$. Writing
$$
\frac{|X-Y|^2}{2n} - 1 \, = \,
\frac{1}{2}\,\Big(\frac{|X|^2}{n} - 1\Big) + 
\frac{1}{2}\,\Big(\frac{|Y|^2}{n} - 1\Big) - 
\frac{1}{n}\left<X,Y\right>,
$$
we also have, by Jensen's inequality,
$$
\Big|\frac{|X-Y|^2}{2n} - 1\Big|^{\frac{3}{2}} \, \leq \,
\Big|\frac{|X|^2 }{n} - 1\Big|^{\frac{3}{2}} + \Big|\frac{|Y|^2 }{n} - 1\Big|^{\frac{3}{2}} +
\frac{2}{n^{3/2}}\, |\left<X,Y\right>|^{\frac{3}{2}}.
$$
Therefore
\bee
\E\,\Big|\frac{|X-Y|^2}{2n} - 1\Big|^{\frac{3}{2}}
 & \leq &
2\,\Big|\frac{|X|^2 }{n} - 1\Big|^{\frac{3}{2}} + 
\frac{2}{n^{3/2}}\, \Big(\E\, |\left<X,Y\right>|^3\Big)^{1/2} \\
 & = &
\frac{2}{n^{3/4}}\,\big(\sigma_3^{3/2} + m_3^{3/2}\big),
\ene
which gives
$$
\bigg(\E\,\Big|\frac{|X-Y|^2}{2n} - 1\Big|^{\frac{3}{2}}\bigg)^{\frac{2}{3}} 
\, \leq \, \frac{2}{\sqrt{n}}\,\big(\sigma_3 + m_3\big).
$$
Hence
\bee
\E\, |\left<X,Y\right>|\, \big|\,e^{-t^2 \big(\frac{|X-Y|^2}{2n} - 1\big)} - 1\big|\,1_B 
 & \leq &
2 t^2 e^{7t^2/8}\,m_3\, (m_3 + \sigma_3) \\
 & \leq &
4 t^2 e^{7t^2/8}\,(m_3^2 + \sigma_3^2),
\ene
and, as a result,
$$
I_2(t) \leq 32\,\frac{e^{-t^2}}{\sqrt{n}}\, (m_3^3 + \sigma_3^3) +
4 t^2 e^{-t^2/8}\,(m_3^2 + \sigma_3^2),
$$
where the factor $e^{-t^2}$ in the first term can be removed without loss of strength. 

Together with the estimate on $I_1(t)$, we get
$$
I(t) \leq \frac{64}{\sqrt{n}}\,(m_3^3 + \sigma_3^3) +
4 t^2 e^{-t^2/8}\,(m_3^2 + \sigma_3^2),
$$
so
$$
\sqrt{I(t)} \, \leq \, \frac{8}{n^{1/4}}\, (m_3^{3/2} + \sigma_3^{3/2}) +
2 |t|\, e^{-t^2/16}\,(m_3 + \sigma_3)
$$
and 
$$
\frac{1}{\sqrt{n}} \int_0^{T_0} \sqrt{I(t)}\,dt  \, \leq \,
\frac{4T_0}{n^{3/4}}\, (m_3^{3/2} + \sigma_3^{3/2}) + 
\frac{C}{\sqrt{n}}\,(m_3 + \sigma_3)
$$
with some absolute constant $C$.

Returning to the bound (7.3), we thus obtain that
$$
c\,\E_\theta\, \rho(F_\theta,F)
 \, \leq \, 
\frac{\log n}{n^{3/4}}\, (m_3^{3/2} + \sigma_3^{3/2}) + 
\frac{C}{\sqrt{n}}\,(m_3 + \sigma_3)
$$
To simplify it, one may use again that $m_3 \geq 1$, which implies that
$m_3 + \sigma_3 \leq 2(m_3^{3/2} + \sigma_3^{3/2})$ for all values of $\sigma_3$.
Thus, with some absolute constant $c>0$,
$$
c\,\E_\theta\, \rho(F_\theta,F) \, \leq \, 
\frac{C}{\sqrt{n}}\,(m_3^{3/2} + \sigma_3^{3/2}).
$$

To get a similar bound with $\Phi$ in place of $F$, i.e. (7.1), one may apply
the estimate $\rho(F,\Phi) \leq C\,\frac{1 + \sigma_2}{\sqrt{n}}$,
where $1 + \sigma_2^2$ may further be bounded by $2(m_3^{3/2} + \sigma_3^{3/2})$.
\qed


\vskip10mm
\section{{\bf The i.i.d. Case}}
\setcounter{equation}{0}

\vskip2mm
\noindent
Theorem 1.3 follows from Theorem 6.3, by taking into account the following
elementary statement (various variants of which under higher moment assumptions 
are well-known).

\vskip5mm
{\bf Lemma 8.1.} {\sl Assume that the non-negative random variables $\xi_1,\dots,\xi_n$
are independent and identically distributed, with $\E \xi_1 = 1$. 
Given $0 < \lambda < 1$, let a number $\kappa > 0$ is chosen to satisfy
$$
\E\,\xi_1\,1_{\{\xi_1 > \kappa\}} \leq \frac{1-\lambda}{2}.
$$
Then for the sum $S_n = \xi_1 + \dots + \xi_n$, we have
\be
\P\{S_n \leq \lambda n\} \leq 
\exp\Big\{- \frac{(1-\lambda)^2}{8 \kappa}\,n\Big\}.
\en
}

{\bf Proof.} Let $V$ denote the common distribution of $\xi_i$. The function 
$$
u(t) = \E\, e^{-t \xi_1} = \int_0^\infty e^{-tx}\,dV(x), \qquad t \geq 0,
$$
is positive, convex, non-increasing, and has a continuous, non-decreasing 
derivative
$$
u'(t) = -\E\, \xi_1 e^{-t \xi_1} = -\int_0^\infty x e^{-tx}\,dV(x),
$$
with $u(0)=1$, $u'(0) = -1$. 

Let $\kappa_p$ denote the maximal quantile for the probability measure $xdV(x)$ 
on $(0,\infty)$ of a given order $p \in (0,1)$, i.e., the minimal number 
such that
$$
\int_{\kappa_p}^\infty x\,dV(x) \leq 1-p,
$$
where the integration is performed over the open half-axis $(\kappa_p,\infty)$.
Using the elementary inequality $1 - e^{-y} \leq y$ ($y \geq 0$), we have, 
for all $s>0$,
\bee
1+u'(s) 
 & = &
\int_0^\infty x (1-e^{-sx})\,dV(x) \\
 & = &
\int_{0 < x \leq \kappa_p} x (1-e^{-sx})\,dV(x) + \int_{x > \kappa_p} 
x (1-e^{-tx})\,dV(x) \\
 & \leq &
s\int_{0 < x \leq \kappa_p} x^2\,dV(x) + p \, \leq \, p + \kappa_p s.
\ene
This gives
\bee
u(t) 
 & = &
1 - t + \int_0^t (1+u'(s))\,ds \\
 & \leq &
1-t + p t + \kappa_p\, \frac{t^2}{2}\, \leq \,
\exp\Big\{-t + p t + \kappa_p\, \frac{t^2}{2}\Big\},
\ene
and therefore
\bee
\P\{S_n \leq \lambda n\} 
 & \leq &
e^{\lambda n t}\,\E\, e^{-t S_n}
 \, = \,
e^{\lambda n t}\, u(t)^n \\ 
 & \leq &
\exp\Big\{-n\Big((1 - \lambda - p) t - 
\kappa_p\, \frac{t^2}{2}\Big)\Big\}.
\ene
If $p < 1-\lambda$, the right-hand side is minimized at
$t = \frac{1 - \lambda - p}{\kappa_p}$, and we get
$$
\P\{S_n \leq \lambda n\} \leq 
e^{-n\,(1 - \lambda - p)^2/2\kappa_p}.
$$
One may take, for example, $p = (1 - \lambda)/2$, and then we arrive at (8.1).
\qed

\vskip5mm
{\bf Proof of Theorem 1.3.} First let us derive the inequality
\be
\E_\theta\, \rho(F_\theta,F) \, \leq \, c\,\sqrt{\frac{\log n}{n}}
\en
with the typical distribution $F$  in place of $G$. Let $Y = (Y_1,\dots,Y_n)$ be 
an independent copy of $X$. Since the Kolmogorov distance is scale invariant, without 
loss of generality one may assume that $\E\,(X_1 - Y_1)^2 = 1$. But then,
by Lemma 8.1, applied to the random variables $\xi_i = (X_i - Y_i)^2$, we have
$$
\P\{|X-Y|^2 \leq n/4\} \leq e^{-cn}
$$
with some constant $c>0$ depending on the distribution of $X_1$ only.
In addition, $M_1 \leq M_2 = \frac{1}{2}$. As a result, Theorem 6.3 yields (8.2).

Now, in order to replace $F$ with $G$ in (8.2), one may apply Proposition 3.1.
Indeed, $F$ represents the distribution function of $rZ_n$, where 
$r = \frac{1}{\sqrt{n}}\,|X|$ and $Z_n = \sqrt{n}\,\theta_1$ is independent of 
$r$. Similarly, $G$ is the distribution function of $rZ$ where $Z \sim N(0,1)$ 
is independent of $r$. Let $\Phi_n$ denote the distribution function of $Z_n$ 
and $\varphi_n$ its density. Since $F(x) = \E\,\Phi_n(x/r)$ and
$G(x) = \E\,\Phi(x/r)$, we conclude, by the triangle inequality, that
\bee
\rho(F,G) 
 & \leq & 
\E\, \sup_x\, |\Phi_n(x/r) - \Phi(x/r)| \\
 & = & 
\sup_x\, |\Phi_n(x) - \Phi(x)| \, \leq \, 
\int_{-\infty}^\infty |\varphi_n(x) - \varphi(x)|\,dx \, \leq \, \frac{C}{n}.
\ene
\qed


\vskip10mm
\begin{thebibliography}{BH3}
\itemsep=-0pt
\small

\vskip2mm 
\bibitem[A-B]{A-B}
C. Aistleitner, and I. Berkes. On the central limit theorem for $f(n_k x)$. 
       Probab. Theory Related Fields  146  (2010),  no. 1-2, 267--289. 

\vskip2mm 
\bibitem[A-E]{A-E}
C. Aistleitner, and C. Elsholtz. The central limit theorem for subsequences 
       in probabilistic number theory. Canad. J. Math. 64 (2012), no. 6, 1201--1221.

\vskip2mm 
\bibitem[Ba]{Ba}
H. Bateman. Higher transcendental functions. Vol. II, McGraw-Hill Book Company, Inc.,
       1953, 396 pp.

\vskip2mm 
\bibitem[B1]{B1}
S. G. Bobkov. On concentration of distributions of random weighted sums. 
       Ann. Probab. 31 (2003), no. 1, 195--215. 
          
\vskip2mm 
\bibitem[B2]{B2}
S. G. Bobkov. On a theorem of V. N. Sudakov on typical distributions. (Russian) 
       J. Math. Sciences (New York), 167 (2010), no. 4, 464--473. Translated from: 
       Zap. Nauchn. Semin. POMI, vol. 368 (2009), 59--74.

\vskip2mm 
\bibitem[B3]{B3}
S. G. Bobkov.  Closeness of probability distributions in terms of Fourier-Stieltjes   
       transforms. Russian Math. Surveys, vol. 71, issue 6 (2016), 1021--1079. 
       Translated from: Uspekhi Matem. Nauk, vol. 71, issue 6 (432), 2016, 37--98.

\vskip2mm 
\bibitem[B-C-G1]{B-C-G1}
S. G. Bobkov, G. P. Chistyakov, and F. G\"otze.  Gaussian mixtures and normal 
       approximation for V. N. Sudakov's typical distributions. Preprint (2017). 
			 To appear in: Zap. Nauchn. Semin. POMI and J. Math. Sciences.

\vskip2mm 
\bibitem[B-C-G2]{B-C-G2}
S. G. Bobkov, G. P. Chistyakov, and F. G\"otze. Concentration and Gaussian 
       approximation for randomized sums. In preparation. 

\vskip2mm 
\bibitem[B-G1]{B-G1}
S. G. Bobkov, and F. G\"otze.  On the central limit theorem along subsequences 
       of noncorrelated observations. Teor. Veroyatnost. i Primenen. 48 (2003), 
			 no. 4, 745--765. Reprinted in: Theory Probab. Appl. 48 (2004), no. 4, 604--621. 

\vskip2mm 
\bibitem[B-G2]{B-G2}
S. G. Bobkov, and F. G\"otze. Concentration inequalities and limit theorems for 
       randomized sums. Probab. Theory Related Fields 137 (2007), no. 1--2, 49--81. 

\vskip2mm
\bibitem[F]{F}
K. Fukuyama. A central limit theorem to trigonometric series with bounded gaps. 
       Probab. Theory Related Fields 149 (2011), no. 1--2, 139--148. 

\vskip2mm
\bibitem[G]{G}
V. F. Gaposhkin. The rate of approximation to the normal law of the distributions 
       of weighted sums of lacunary series. (Russian) 
       Teor. Verojatnost. i Primenen. 13 (1968), 445--461. 

\vskip2mm
\bibitem[Ka]{Ka}
M. Kac. On the distribution of values of sums of the type $\sum f(2^k t)$.
       Ann. Math. (2) 47 (1946), 33--49. 

\vskip2mm 
\bibitem[Kl]{Kl}
B. Klartag. A Berry-Esseen type inequality for convex bodies with an 
       unconditional basis. Probab. Theory Related Fields 145 (2009), no. 1-2, 
			 1--33. 

\vskip2mm 
\bibitem[K-S]{K-S}
B. Klartag, and S. Sodin. Variations on the Berry-Esseen theorem. Teor. Veroyatn. 
       Primen. 56 (2011), no. 3, 514--533;  reprinted in: Theory Probab. Appl.  
			 56 (2012), no. 3, 403--419.

\vskip2mm 
\bibitem[L]{L}
 M. Ledoux. Concentration of measure and logarithmic Sobolev inequalities.
       S\'eminaire de Probabilit\'es XXXIII. Lect. Notes in Math.
       1709 (1999), 120--216, Springer.

\vskip2mm 
\bibitem[M]{M} 
V. K. Matskyavichyus. A lower bound for the rate of convergence in the central 
       limit theorem. (Russian) Teor. Veroyatnost. i Primenen. 28 (1983), 
			 no. 3, 565--569.

\vskip2mm 
\bibitem[M-M]{M-M} 
E. S. Meckes, and M. W. Meckes. The central limit problem for random vectors 
       with symmetries. J. Theoret. Probab. 20 (2007), no. 4, 697--720. 

\vskip2mm 
\bibitem[P1]{P1} 
V. V. Petrov. Sums of independent random variables. Translated from the Russian 
       by A. A. Brown. Ergebnisse der Mathematik und ihrer Grenzgebiete, Band 82. 
       Springer--Verlag, New York--Heidelberg, 1975. x+346 pp.

\vskip2mm 
\bibitem[P2]{P2} 
V. V. Petrov. Limit theorems for sums of independent random variables (Russian), 
       Nauka, Moscow, 1987. 318 pp.          
         
\vskip2mm
\bibitem[S-Z1]{S-Z1}
R. Salem, and A. Zygmund. On lacunary trigonometric systems.
       Proc. Nat. Acad. Sci. USA, 33 (1947), 333--338.

\vskip2mm
\bibitem[S-Z2]{S-Z2}
R. Salem, and A. Zygmund. On lacunary trigonometric series. II. 
       Proc. Nat. Acad. Sci. U. S. A. 34 (1948), 54--62.

\vskip2mm
\bibitem[S]{S}   
V. N. Sudakov. Typical distributions of linear functionals in finite-dimensional 
       spaces of high dimension. (Russian) Soviet Math. Dokl. 19 (1978), 
			 1578--1582; translation in: Dokl. Akad. Nauk SSSR, 243 (1978), 
			 no. 6, 1402--1405.        


\end{thebibliography}
\end{document}